\documentclass[11pt]{article}
\usepackage{amsfonts}
\usepackage{amssymb}

\newtheorem{theorem}{Theorem}[section]

\newtheorem{lemma}[theorem]{Lemma}

 \usepackage{amssymb}
\usepackage{epsfig}
\usepackage{amsmath}
\usepackage{amsfonts}

\newcommand{\R}{{{\Bbb R}}}

\newcommand{\N}{{{\Bbb N}}}
\newcommand{\Nn}{{\scriptsize {\Bbb N}}}

\def\qed{\hbox to 0pt{}\hfill$\rlap{$\sqcap$}\sqcup$\medbreak}

\title{A simple proof of the Fundamental Theorem of Calculus for the Lebesgue integral}
\date{March, 2012}
\begin{document}
\maketitle

\vspace{1cm}

\begin{center}
{\large Rodrigo L\'opez Pouso \\
Departamento de
An\'alise Matem\'atica\\
Facultade de Matem\'aticas,\\Universidade de Santiago de Compostela, Campus Sur\\
15782 Santiago de
Compostela, Spain.
%\\
%Phone: 34 981 56 31 00,
%Ext. 13213 \\ FAX: 34 981 59 70 54
}
\end{center}

\begin{abstract}
This paper contains a new elementary proof of the Fundamental Theorem of Calculus for the Lebesgue integral. The hardest part of our proof simply concerns the convergence in ${\rm L}^1$ of a certain sequence of step functions, and we prove it using only basic elements from Lebesgue integration theory. 
 \end{abstract}

 \section{Introduction}
 Let $f:[a,b]\longrightarrow \R$ be absolutely continuous on $[a,b]$, i.e., for every $\varepsilon>0$ there exists $\delta>0$ such that if $\{(a_j,b_j)\}_{j=1}^{n}$ is a family of pairwise disjoint subintervals of $[a,b]$ satisfying $$\sum_{j=1}^{n}{(b_j-a_j)}<\delta$$
then
$$\sum_{j=1}^{n}{|f(b_j)-f(a_j)|}<\varepsilon.$$

Classical results ensure that $f$ has a finite derivative almost everywhere in $I=[a,b]$, and that $f'\in {\rm L}^1(I)$, see \cite{bot} or \cite[Corollary 6.83]{str}. These results, which we shall use in this paper, are the first steps in the proof of the main connection between absolute continuity and Lebesgue integration: the Fundamental Theorem of Calculus for the Lebesgue integral. 

\begin{theorem}
\label{tfc}
If $f:I=[a,b] \longrightarrow \R$ is absolutely continuous on $I$ then
$$f(b)-f(a)=\int_{a}^b{f'(x) \, dx} \quad \mbox{in Lebesgue's sense.}$$
\end{theorem}

In this note we present a new elementary proof to Theorem \ref{tfc} which seems more natural and easy than the existing ones. Indeed, our proof can be sketched simply as follows:
\begin{enumerate}
\item We consider a well--known sequence of step functions $\{h_n\}_{n \in \Nn}$ which tends to $f'$ almost everywhere in $I$ and, moreover,
$$\int_a^b{h_n(x) \, dx}=f(b)-f(a) \quad \mbox{for all $n \in \N$.}$$
\item We prove, by means of elementary arguments, that
$$\lim_{n \to \infty}\int_a^b{h_n(x) \, dx}=\int_a^b{f'(x) \, dx}.$$
\end{enumerate}

More precise comparison with the literature on Theorem \ref{tfc} and its several proofs will be given in Section 3.

In the sequel $m$ stands for the Lebesgue measure in $\R$.

\section{Proof of Theorem \ref{tfc}}
For each $n \in \N$ we consider the partition of the interval $I=[a,b]$ which divides it into $2^n$ subintervals of length $(b-a)2^{-n}$, namely
$$x_{n,0}<x_{n,1}<x_{n,2}<\cdots<x_{n,2^n},$$
where $x_{n,i}=a+i(b-a)2^{-n}$ for $i=0,1,2,\dots, 2^n$.

Now we construct a step function $h_n:[a,b) \longrightarrow \R$ as follows: for each $x \in [a,b)$ there is a unique $i \in \{0,1,2,\dots, 2^n-1\}$ such that
$$x \in [x_{n,i},x_{n,i+1}),$$
and we define
$$h_n(x)=\dfrac{f(x_{n,i+1})-f(x_{n,i})}{x_{n,i+1}-x_{n,i}}=\dfrac{2^n}{b-a}[f(x_{n,i+1})-f(x_{n,i})].$$

On the one hand, the construction of $\{h_n\}_{n \in \Nn}$ implies that
\begin{equation}
\label{derlim}
\lim_{n \to \infty}{h_n(x)}=f'(x) \quad \mbox{for all $x \in [a,b) \setminus N$,}
\end{equation}
where $N \subset I$ is a null--measure set such that $f'(x)$ exists for all $x \in I \setminus N$. 

On the other hand, for each $n \in \N$ we compute
\begin{equation}
\nonumber
\int_a^b{h_n(x)  dx}=\sum_{i=0}^{2^n-1}{\int_{x_{n,i}}^{x_{n,i+1}}{h_n(x) dx}}=\sum_{i=0}^{2^n-1}[f(x_{n,i+1})-f(x_{n,i})]=f(b)-f(a),
\end{equation}
and therefore it only remains to prove that
$$\lim_{n \to \infty}\int_a^b{h_n(x) \, dx}=\int_a^b{f'(x) \, dx}.$$
Let us prove that, in fact, we have convergence in ${\rm  L}^1(I)$, i.e.,
\begin{equation}
\label{goal}
\lim_{n \to \infty}\int_a^b{|h_n(x)-f'(x)| \, dx}=0.
\end{equation}
 
Let $\varepsilon>0$ be fixed and let $\delta>0$ be one of the values corresponding to $\varepsilon/4$ in the definition of absolute continuity of $f$. 

Since $f' \in {\rm L}^1(I)$ we can find $\rho>0$ such that for any measurable set $E \subset I$ we have
\begin{equation}
\label{eq1}
\int_E{|f'(x)| \, dx}<\dfrac{\varepsilon}{4} \quad \mbox{whenever $m(E)<\rho$.}
\end{equation}

The following lemma will give us fine estimates for the integrals when $|h_n|$ is ``small". We postpone its proof for better readability.

\begin{lemma}
\label{le1}
For each $\varepsilon>0$ there exist $k, n_k \in \N$ such that 
$$k \cdot m \left( \left\{x \in I \, : \, \sup_{n \ge n_k}|h_n(x)|>k \right\} \right)<\varepsilon.$$
\end{lemma}

Lemma \ref{le1} guarantees that there exist $k, n_k \in \N$ such that 
\begin{equation}
\label{etdos}
k \cdot m \left( \left\{x \in I \, : \, \sup_{n \ge n_k}|h_n(x)|>k \right\} \right)<\min\left\{\delta,\dfrac{\varepsilon}{4}, \rho \right\}.
\end{equation}
Let us denote
$$A=\left\{x \in I \, : \, \sup_{n \ge n_k}|h_n(x)|>k \right\},$$
which, by virtue of (\ref{etdos}) and (\ref{eq1}), satisfies the following properties:
\begin{eqnarray}
\label{intdos}
m(A) <\delta,   \\
\label{et2}
k\cdot m(A)<\dfrac{\varepsilon}{4},\\
\label{int1}
\int_{A}{|f'(x)| \, dx}<\dfrac{\varepsilon}{4}.
\end{eqnarray}

We are now in a position to prove that the integrals in (\ref{goal}) are smaller than $\varepsilon$ for all sufficiently large values of $n \in \N$. We start by noticing that (\ref{int1}) guarantees that for all $n \in \N$ we have  
\begin{align}
\nonumber \int_{I}{|h_n(x)-f'(x)| \, dx}&=\int_{I \setminus A}{|h_n(x)-f'(x)| \, dx}+\int_{ A}{|h_n(x)-f'(x)| \, dx}\\
\label{int2}
&<\int_{I \setminus A}{|h_n(x)-f'(x)| \, dx}+\int_{ A}{|h_n(x)| \, dx}+\dfrac{\varepsilon}{4}.
\end{align}

The definition of the set $A$ implies that for all $n \in \N$, $n \ge n_k$, we have
$$|h_n(x)-f'(x)| \le k+|f'(x)| \quad \mbox{for almost all $x \in I \setminus A$},$$
so the Dominated Convergence Theorem yields
\begin{equation}
\label{et3}
\lim_{n \to \infty}\int_{I \setminus A}{|h_n(x)-f'(x)| \, dx}=0.
\end{equation}

From (\ref{int2}) and (\ref{et3}) we deduce that there exists $n_{\varepsilon} \in \N$, $n_{\varepsilon} \ge n_k$, such that for all $n \in \N$, $n \ge n_{\varepsilon}$, we have
\begin{equation}
\label{int3}
\int_{I}{|h_n(x)-f'(x)| \, dx}<\dfrac{\varepsilon}{2}+\int_{ A}{|h_n(x)| \, dx}.
\end{equation}

Finally, we estimate $\int_A{|h_n|}$ for each fixed $n \in \N$, $n \ge n_{\varepsilon}$. First, we decompose $A=B \cup C$, where
$$B=\{ x \in A \, : \, |h_n(x)| \le k\} \quad \mbox{and} \quad C=A \setminus B.$$

We immediately have
\begin{equation}
\label{int4}
\int_B{|h_n(x)| \, dx} \le k \cdot m(B) \le k \cdot m(A) < \dfrac{\varepsilon}{4} \quad \mbox{by (\ref{et2}).}
\end{equation}

Obviously, $\int_C{|h_n|}<\varepsilon/4$ when $C=�\varnothing$. Let us see that this inequality holds true when $C \neq \varnothing$. For every $x \in C=\{ x \in A \, : \, |h_n(x)| >k \}$ there is a unique index $i \in \{0,1,2,\dots, 2^n-1\}$ such that $x \in [x_{n,i},x_{n,i+1})$. Since $|h_n|$ is constant on  $[x_{n,i},x_{n,i+1})$ we deduce that  $[x_{n,i},x_{n,i+1}) \subset C$. Thus there exist indexes $i_l \in \{0,1,2,\dots, 2^n-1\}$, with $l=1,2,\dots, p$ and $i_l \ne i_{\tilde l}$ if $l \neq \tilde l$, such that
$$C=\bigcup_{l=1}^{p}[x_{n,i_l},x_{n,i_l+1}).$$
Therefore
$$\sum_{l=1}^{p}(x_{n,i_l+1}-x_{n,i_l})=m(C) \le m(A) < \delta \quad \mbox{by (\ref{intdos}),}$$
and then the absolute continuity of $f$ finally comes into action:
\begin{align*}
\int_C{|h_n(x)| \, dx}&=\sum_{l=1}^{p}\int_{x_{n,i_l}}^{x_{n,i_l+1}}{|h_n(x)| \, dx} \\
&=\sum_{l=1}^{p}|f(x_{n,i_l+1})-f(x_{n,i_l})|< \dfrac{\varepsilon}{4}.
\end{align*}

This inequality, along with (\ref{int3}) and (\ref{int4}), guarantee that for all $n \in \N$, $n \ge n_{\varepsilon}$, we have
$$\int_{I}{|h_n(x)-f'(x)| \, dx}<\varepsilon,$$
thus proving (\ref{goal}) because $\varepsilon$ was arbitrary. \qed

\noindent
{\bf Proof of Lemma \ref{le1}.} Let $\varepsilon>0$ be fixed and let $\rho>0$ be such that for every measurable set $E \subset I$ with $m(E)< \rho$ we have
$$\int_{E}{|f'(x)| \, dx}< \dfrac{\varepsilon}{2}.$$
Let $N \subset I$ be as in (\ref{derlim}) and let $k \in \N$ be sufficiently large so that
$$m \left( \left\{x \in I \setminus N \, : \, |f'(x)| \ge k  \right\}\right) < \rho,$$
which implies that
\begin{equation}
\label{l1}
k \cdot m \left( \left\{x \in I \setminus N \, : \, |f'(x)| \ge k  \right\}\right) \le \int_{\left\{x \in I \setminus N \, : \, |f'(x)| \ge k  \right\}}{|f'(x)| \, dx}<\dfrac{\varepsilon}{2}.
\end{equation}
Let us define
$$E_j=\left\{ x \in I \setminus N \, : \,  \sup_{n \ge j}|h_n(x)| >k \right\} \quad (j \in \N).$$
Notice that $E_{j+1}\subset E_j$ for every $j \in \N$, and $m(E_1)<\infty$, hence
\begin{equation}
\label{pl1}
\lim_{j \to \infty}m(E_j)=m \left( \bigcap_{j=1}^{\infty}E_j\right).
\end{equation}

Clearly, $\cap_{j=1}^{\infty}E_j \subset \left\{x \in I \setminus N \, : \, |f'(x)| \ge k  \right\}$, so we deduce from (\ref{pl1}) that we can find some $n_k \in \N$ such that  
$$m(E_{n_k}) \le m \left( \left\{x \in I \setminus N \, : \, |f'(x)| \ge k  \right\}\right) +\dfrac{\varepsilon}{2k},$$
and then (\ref{l1}) yields $k \cdot m(E_{n_k}) < \varepsilon$. \qed
\section{Final remarks}
The sequence $\{h_n\}_{n \in \Nn}$ is used in other proofs of Theorem \ref{tfc}, see \cite{bar} or \cite{vol}. The novelty in this paper is our elementary and self--contained proof of (\ref{goal}). Incidentally, a revision of the proof of our Lemma \ref{le1} shows that it holds true for any sequence of measurable functions $\tilde h_n: E \subset \R \longrightarrow \R$ which converges pointwise almost everywhere to some $h \in {\rm L}^1(E)$ and $m(E) < \infty$.
 
Our proof avoids somewhat technical results often invoked to prove Theorem \ref{tfc}. For instance, we do not use any sophisticated estimate for the measure of image sets such as \cite[Theorem 7.20]{bbt}, \cite[Lemma 6.88]{str} or \cite[Proposition 1.2]{vol}, see also \cite{kol}. We do not use the following standard lemma either: an absolutely continuous function having zero derivative almost everywhere is constant, see \cite[Theorem 7.16]{bbt} or \cite[Lemma 6.89]{str}. It is worth having a look at \cite{gor} for a proof of that lemma using tagged partitions; see also \cite{bot1} for a proof  based on full covers \cite{tho}. Concise proofs of Theorem \ref{tfc} follow from the Radon--Nikodym Theorem, see \cite{bar}, \cite{bbt} or \cite{rud}, but this is far from being elementary.   

Finally, it is interesting to note that (\ref{goal}) easily follows from the Dominated Convergence Theorem when $f$ is Lipschitz continuous on $I$. This fact made the author think about the following project for students in an introductory course to Lebesgue integration.

\begin{center}
{\sc Project: Two important results for the price of one.}
\end{center}

\begin{enumerate}
\item Let $f:I=[a,b] \longrightarrow \R$ be Lispchitz continuous on $I$. A deep result (worth to know without proof) guarantees that $f'(x)$ exists for almost all $x \in I$, see \cite[Theorem 7.8]{bbt}.
 
 Consider the sequence $\{h_n\}_{n \in \Nn}$ as defined in Section 2 and prove 
\begin{enumerate}
\item $\{h_n (x) \}_{n \in \Nn}$ tends to $f'(x)$ for almost all $x \in I$;
\item $\int_I{h_n(x) \, dx}=f(b)-f(a)$ for all $n \in \N$;
\item (Use the Dominated Convergence Theorem) $f' \in L^1(I)$ and 
$$f(b)-f(a)=\int_a^b{f'(x) \, dx} \quad \mbox{in Lebesgue's sense.}$$
\end{enumerate}
\item Let $g:I=[a,b] \longrightarrow \R$ be Riemann--integrable on $I$ and define
$$f(x)=(R)\int_a^x{g(s) \, ds} \quad (x \in I),$$
where $(R) \int$ stands for the Riemann integral.

Use the information in Exercise 1 to deduce that $g \in {\rm L}^1(I)$ and
$$(R) \int_a^b{g(x) \, dx}=\int_a^b{g(x) \, dx} \quad \mbox{in Lebesgue's sense.}$$

\end{enumerate}

\end{document}